\documentclass{article}
\usepackage{amssymb,amsmath}
\usepackage{graphicx,CJK}

\def\qed{\hfill {\hbox{${\vcenter{\vbox{               
   \hrule height 0.4pt\hbox{\vrule width 0.4pt height 6pt
   \kern5pt\vrule width 0.4pt}\hrule height 0.4pt}}}$}}}

\def\bar{\overline}
\def\trace{\mathrm{tr\ }}

\newtheorem{theorem}{Theorem}
\newtheorem{definition}{Definition}

\newtheorem{example}{Example}

\newtheorem{remark}{Remark}

{\qed\par\medskip}

\date{}

\title{\Large \textbf{Quantum Enhancements of Involutory Birack Counting 
Invariants}}

\author{Sam Nelson\footnote{Email: knots@esotericka.org}\and
Veronica Rivera\footnote{Email: verorive1@gmail.com}}

\begin{document}
\maketitle

\begin{abstract}
The involutory birack counting invariant is an integer-valued invariant of
unoriented tangles defined by counting homomorphisms from the fundamental
involutory birack of the tangle to a finite involutory birack over a set 
of framings modulo the birack rank of the labeling birack. In this first
of an anticipated series of several papers, we enhance the involutory birack 
counting invariant with \textit{quantum weights}, which may be understood as 
tangle functors of involutory birack-labeled unoriented tangles. 
\end{abstract}

\parbox{5.5in} {\textsc{Keywords:} Biracks, enhancements of counting invariants,
quantum invariants
\smallskip

}

\section{\large\textbf{Introduction}}\label{I}

A \textit{birack} $\mathbf{X}$ is an algebraic structure with axioms determined 
by the framed Reidemeister moves such that labelings of arcs in a 
tangle diagram with elements of $\mathbf{X}$ before and after a move are in 
bijective correspondence. In particular, the number of $\mathbf{X}$-labelings 
of a tangle diagram by a given involutory birack $\mathbf{X}$ is an invariant
of unoriented framed links under framed isotopy moves. Biracks were introduced 
in \cite{FRS}; 
special cases of biracks including racks, quandles and biquandles have been 
studied and used to define link invariants in various works including 
\cite{FR,J,M,CJKLS,KR,FJK} and more. 

A finite birack has an associated integer $N$ called the \textit{birack rank} 
or \textit{birack characteristic} with the property that framed unoriented 
links with equivalent framings modulo $N$ have equal numbers of labelings.
Summing over a complete period of framings, one obtains an integer-valued
invariant of unframed tangles called the \textit{integral birack counting 
invariant}, denoted $\Phi_{\mathbf{X}}^{\mathbb{Z}}$; see \cite{N} for more.

An \textit{enhancement} of $\Phi_{\mathbf{X}}^{\mathbb{Z}}$ is obtained by 
evaluating an invariant $\sigma$ of $\mathbf{X}$-labeled diagrams on each 
$\mathbf{X}$-labeling of a tangle diagram $T$ and collecting these values 
over a complete set of $\mathbf{X}$-labelings to obtain a multiset of 
$\sigma$-values sometimes called ``signatures''. This multiset is then a
generally stronger invariant of unframed links whose cardinality is the
birack counting invariant $\Phi_{\mathbf{X}}^{\mathbb{Z}}$.

In this paper we develop an enhancement of $\Phi_{\mathbf{X}}^{\mathbb{Z}}$ in the
special case when $\mathbf{X}$ is an \textit{involutory birack}, the type of 
birack appropriate for defining invariants of unoriented tangles, using 
$\sigma$ values we call \textit{quantum weights}. Quantum weights are 
birack-labeled tangle functors which may be understood as customized 
quantum invariants for $\mathbf{X}$-labeled tangle diagrams. 
This will be the first is a series of papers in which we consider progressively
more complex cases, starting with the unoriented case in this paper and in 
the sequels considering the oriented and virtual cases. The paper is organized 
as follows. In Section \ref{IR} we review the basics of involutory biracks 
and the birack counting invariant. In Section \ref{QE} we introduce quantum 
weights and the quantum-enhanced birack counting invariant; we give examples 
to demonstrate how the invariant is computed and show
that the enhanced invariant is stronger than the unenhanced invariant.
In Section \ref{B} we consider quantum enahancements of involutory birack 
counting invariants for closed braids.
We conclude in Section \ref{Q} with some open questions and directions for 
future work.

\section{\large\textbf{Involutory Biracks and the Counting Invariant}}\label{IR}

Recall that a \textit{framed unoriented tangle}  is an equivalence class of 
disjoint unions of simple closed curves and arcs in 
$A=\mathbb{R}^2\times [0,1]$ with endpoints in fixed positions in $\partial A$;
two such disjoint unions are equivalent if they can be connected by an ambient 
isotopy of $A$ fixing $\partial A$ and preserving the linking number of each 
component with a 
choice of framing curve for each component. The linking number of a component
with its framing curve is the \textit{framing number} of the component; the 
framing number is equal to the \textit{writhe} or sum of crossing signs 
\[\includegraphics{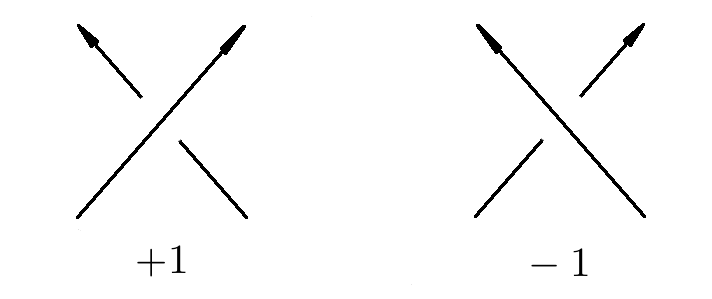}\]
at self-crossings when the framing curve is the \textit{blackboard framing 
curve} obtained by pushing off a parallel copy of each component in a diagram
of $T$. 

Equivalently, framed unoriented tangles can be defined combinatorially
as an equivalence class of unoriented tangle diagrams under the equivalence
relation generated by the \textit{framed unoriented Reidemeister moves}
\[
\includegraphics{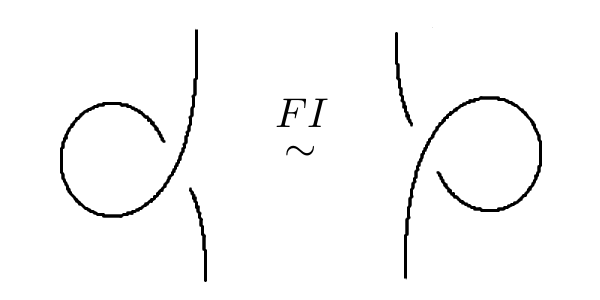} \quad
\includegraphics{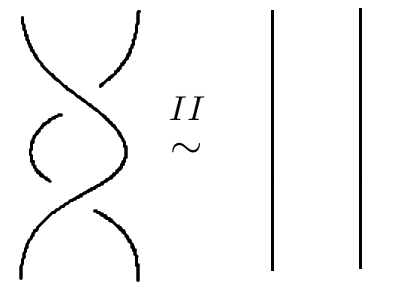} \quad \includegraphics{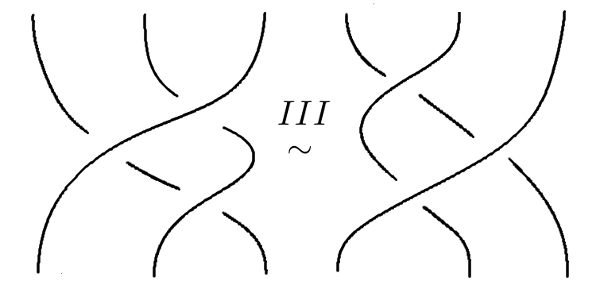}
\]
In particular, the usual unframed Reidemeister I move changes the backboard 
framing of a component by $\pm 1$, while the framed version depicted above 
preserves the framing.
An unframed tangle of $c$ components thus determines a $\mathbb{Z}^c$--lattice 
of framed tangles; a choice of ordering on the components allows us to specify
a framing with a \textit{framing vector} $\vec{w}\in \mathbb{Z}^c$.

Next, we have a definition from \cite{AN}:

\begin{definition}
\textup{An \textit{involutory birack} is a set $\mathbf{X}$ with an 
invertible map 
$B:\mathbf{X}\times \mathbf{X}\to \mathbf{X}\times \mathbf{X}$ such that
\begin{list}{}{}
\item[(i)] $(\tau B)^2=I$
\item[(ii)] The components $(\tau B\Delta)_{1,2}:\mathbf{X}\to\mathbf{X}$ 
of the map $\tau B\Delta:\mathbf{X}\to \mathbf{X}\times\mathbf{X}$ are 
bijections, and
\item[(iii)] B satisfies the \textit{set-theoretic Yang-Baxter equation}
\[(B\times I)(I\times B)(B\times I)=(I\times B)(B\times I)(I\times B)\]
\end{list}
where $\tau:\mathbf{X}\times \mathbf{X}\to \mathbf{X}\times \mathbf{X}$, 
$\Delta:\mathbf{X}\to \mathbf{X}\times \mathbf{X}$ and 
$I:\mathbf{X}\to \mathbf{X}$ are defined by $\tau(x,y)=(y,x)$, 
$\Delta(x)=(x,x)$ and $I(x)=x$ respectively. The map 
$S:\mathbf{X}\times \mathbf{X}\to \mathbf{X}\times \mathbf{X}$ defined by 
$S=\tau B\tau= B^{-1}$ is called the \textit{sideways map}. We will find it 
useful to abbreviate $B_1(x,y)=y^x$, $B_2(x,y)=x_y$, $(S\Delta)_2^{-1}=\alpha$ 
and $(S\Delta)_1\alpha=\pi$.}
\end{definition}

Let $\mathbf{X}$ be an involutory birack.  A \textit{birack labeling} or 
\textit{$\mathbf{X}$--labeling}
of an unoriented framed tangle diagram $T$ is an assignment of elements of
$\mathbf{X}$ to the \textit{semiarcs} of $T$, i.e., the portions of $T$ between
crossing points, such that at every crossing we have the pictured relationship
between semiarc labelings:
\[\includegraphics{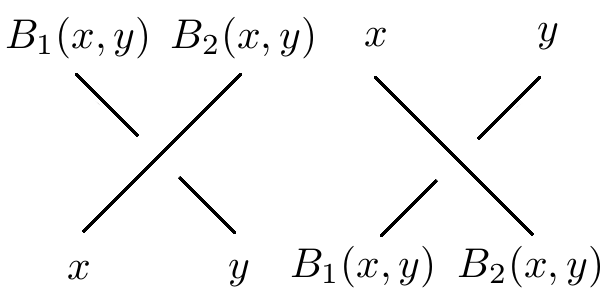}.\]
The involutory birack axioms are then consequences of the framed Reidemeister 
moves; more precisely, they are the conditions required to guarantee that for 
every birack labeling of a tangle diagram before a framed unoriented 
Reidemeister move, there is exactly one corresponding labeling of the diagram 
$T'$ after the Reidemeister move. In particular, given a finite involutory 
birack $\mathbf{X}$, the number of $\mathbf{X}$--labelings of a framed tangle 
$T$ is an invariant of unoriented framed isotopy called the \textit{basic 
birack counting invariant}, denoted $\Phi_{\mathbf{X}}^B(L)$. 

Standard examples of involutory biracks include
\begin{itemize}
\item \textit{Constant Action Involutory Biracks}. Let $\mathbf{X}$ be any 
set and $\sigma,\tau:\mathbf{X}\to \mathbf{X}$ two involutions such that
$\sigma\tau=\tau\sigma$. Then \[B(x,y)=(\sigma(y),\tau(x))\] 
defines an involutory birack structure on $\mathbf{X}$.
\item \textit{Involutory $(t,s,r)$-Biracks.} Let $\mathbf{X}$ be a module 
over the ring
$\tilde\Lambda=\mathbb{Z}[t,s,r]/I$ where $I$ is the ideal generated by 
$s^2-s(1-tr),1-t^2,1-r^2,(t+r)s,$ and $(1-r)s$. Then $\mathbf{X}$ is an 
involutory birack with map
\[B(x,y) = (ty+sx,rx).\]
\item \textit{The Fundamental Involutory Birack of an Unoriented Framed Tangle.}
Let $T$ be an unoriented framed tangle and let $G=\{g_1,\dots,g_n\}$ be a set 
of generators corresponding bijectively with the semiarcs in $T$. The set 
$W$ of \textit{involutory birack words} in $T$ is defined recursively by the
rules (1) $x\in G\Rightarrow x\in W$ and (2) 
$x,y\in W\Rightarrow B_1(x,y)\in W$ and $B_2(x,y)\in W$. Then the 
\textit{Fundamental Involutory Birack} of $T$, $IB(T)$, is the set of 
equivalence classes in $W$ under the equivalence relation determined by
the involutory birack axioms and the crossing relations in $T$. Note that 
birack labelings of a tangle $T$ by a birack $\mathbf{X}$ are precsiely birack
homomorphisms $f:IB(T)\to \mathbf{X}$, i.e. maps satisfying 
$(f\times f)B_{IB(T)}=B_{\mathbf{X}}(f\times f)$ where $B_{\mathbf{X}}$ and $B_{IB(T)}$
are the birack maps in $\mathbf{X}$ and $IB(T)$ respectively.
\end{itemize}

Given a finite set $\mathbf{X}=\{x_1,\dots,x_n\}$ we can define an involutory
birack structure on $\mathbf{X}$ by giving a matrix $M_{\mathbf{X}}=[U|L]$ 
encoding the operation tables of the components of $B$ considered as binary 
operations $B(x,y)=(y^x,x_y).$ That is, $M_{\mathbf{X}}$ is a block matrix 
with two blocks $U$ and $L$ such that the entries in row $i$, column $j$ of 
$U^T$ and $L$ respectively are $k$ and $l$ where 
$B(x_i, x_j)=(x_k,x_l)$.\footnote{We use the transpose of $U$ so that in an 
$\mathbf{X}$-labeling of a tangle diagram, the row label operand and the 
output label lie on the same strand.} Such a matrix defines an involutory 
birack if and only if the map $B:\mathbf{X}\times \mathbf{X}\to 
\mathbf{X}\times \mathbf{X}$ defined by the matrix satisfies the 
birack axioms.

\begin{example}\label{ex0}\textup{
The smallest involutory birack which is neither a rack nor a quandle is 
given by the birack matrix
\[\left[\begin{array}{cc|cc}
1 & 1 & 2 & 2 \\
2 & 2 & 1 & 1
\end{array}
\right].\]
}\end{example}

For a birack $\mathbf{X}$, the \textit{kink map} $\pi:\mathbf{X}\to \mathbf{X}$
defined by $\pi=(S\Delta)_1(S\Delta)_2^{-1}$ is a bijection representing 
going through a positive kink 
as pictured.
\[\includegraphics{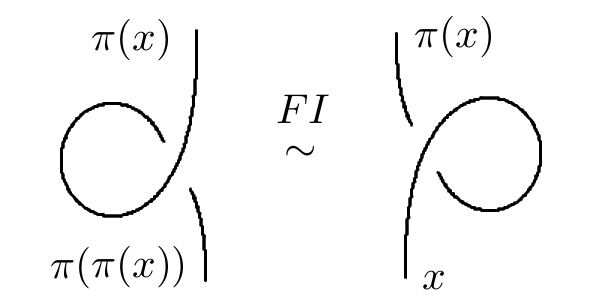}\]
The exponent of $\pi$ considered as an element of the symmetric group 
$S_{\mathbf{X}}$, i.e. the smallest positive integer $N$ such that 
$\pi^N=I$, is called the \textit{birack rank} or \textit{birack characteristic} of $\mathbf{X}$. In the unoriented case, the framed type I move then requires 
that $\pi=\pi^{-1}$ and we obtain
\begin{theorem}
An involutory birack has birack rank $N=1$ or $N=2$.
\end{theorem}
An involutory birack of rank $N=1$ is an \textit{involutory biquandle} or 
\textit{bikei} (\begin{CJK*}{UTF8}{min}双圭\end{CJK*}).

By construction, if $\mathbf{X}$ is an involutory birack of rank $N$, then 
$\mathbf{X}$-labelings of a framed tangle diagram before and after the 
\textit{$N$-phone cord move}
are in bijective correspondence. When $N=1$, the $N$-phone cord move is the
unframed Reidemeister type I move; when $N=2$, the move is as pictured below.
\[\includegraphics{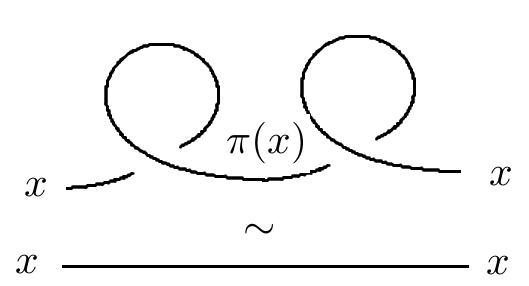}\]

Given an unoriented framed tangle $T$ of $c$ components and a finite involutory
birack $\mathbf{X}$, the $\mathbb{Z}^c$--lattice of unoriented framings of $T$ 
determines a 
$\mathbb{Z}^c$--lattice of basic counting invariant values 
$\Phi_{\mathbf{X}}^B(T,\vec{w})$. 
\[\includegraphics{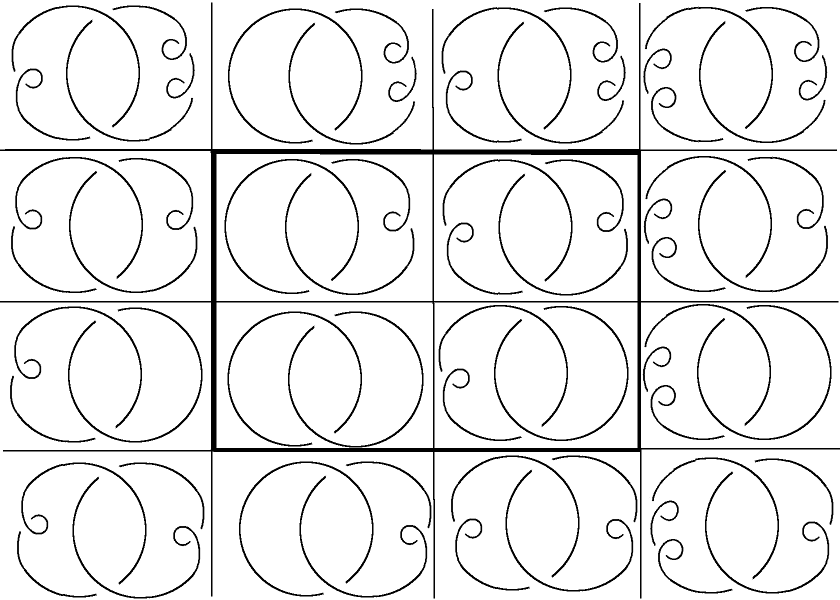}\quad 
\raisebox{0.35in}{\includegraphics{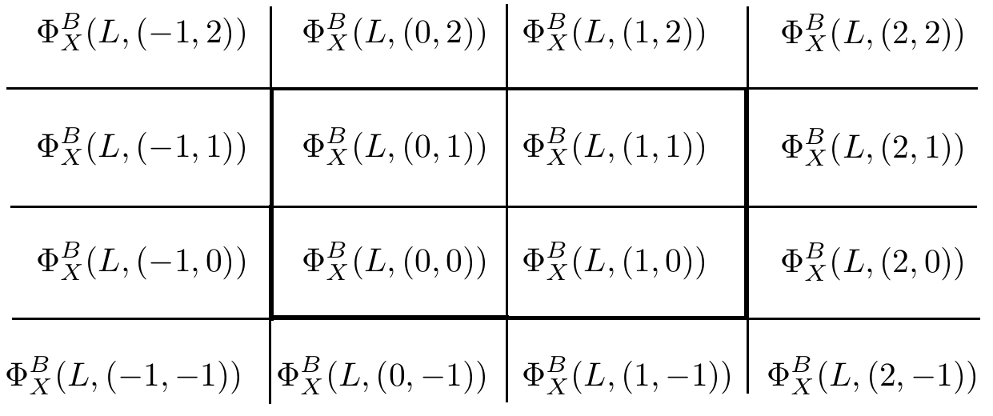}}
\]
If $N=1$, then these basic counting invariants are all equal; if $N=2$, then 
any two framings of
$T$ with framing vectors congruent mod 2 have the same basic counting 
invariants, so the $\mathbb{Z}^c$--lattice is tiled with a $2\times 2$ tile
of basic counting invariant values. In either case, we can sum the basic 
counting invariants over an $N\times N$ tile of framings to obtain an invariant
of unframed unoriented links $T$ called the \textit{integral involutory birack
counting invariant}, 
\[\Phi^{\mathbb{Z}}_{\mathbf{X}}(T)
=\sum_{\vec{w}\in(\mathbb{Z}_N)^C}\Phi_{\mathbf{X}}^{B}(T,\vec{w})\]
where $(T,\vec{w})$ is a diagram of $T$ with framing vector $\vec{w}$.

\begin{example}\textup{The matrix 
\[M_{\mathbf{X}}=\left[\begin{array}{rrr|rrr}
2 & 2 & 2 & 1 & 1 & 1 \\
1 & 1 & 1 & 2 & 2 & 2 \\
3 & 3 & 3 & 3 & 3 & 3 \\
\end{array}\right]\]
defines an involutory birack of rank $N=2$. To compute the counting invariant 
$\Phi_{\mathbf{X}}^{\mathbb{Z}}(T)$ of a tangle $T$,
we need to find the rack labelings over a complete set of framing
vectors mod $2$. For example, the Hopf link $L2a1$ has two components
and thus the space of framings is 
$(\mathbb{Z}_n)2^2=\{(0,0),(1,0),(0,1),(1,1)\}$. 
The Hopf link with framing vector $(0,1)$, for instance, has 
$\mathbf{X}$-labelings
\[\includegraphics{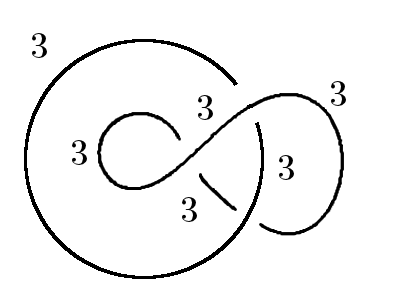}
\includegraphics{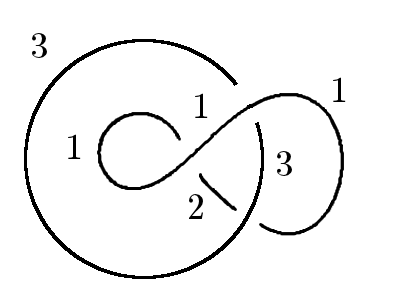}
\includegraphics{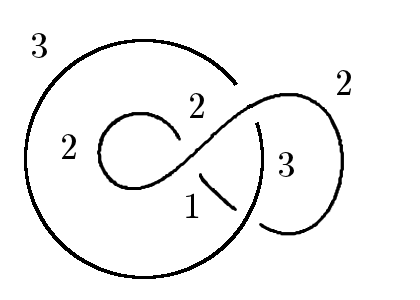}
\]
so this framing contributes $3$ labelings to the invariant.
The other framing vectors contribute $1,3,$ and $5$ labelings respectively,
so we have $\phi_{\mathbf{X}}^{\mathbb{Z}}(L2a1)= 1+3+3+5=12$.
}\end{example}

\section{\large\textbf{Quantum Enhancements}}\label{QE}

Let $T$ be an unoriented link of $c$ components and $\mathbf{X}$ a finite 
involutory birack of rank $N$. To each $\mathbf{X}$-labeling 
$f:IB(T)\to \mathbf{X}$ of $T$ we would like to define a signature 
$\sigma(f)$ which is invariant under $\mathbf{X}$-labeled
framed Reidemeister moves. We will do this by defining an 
$\mathbf{X}$-labeled tangle functor or $\mathbf{X}$-labeled 
quantum invariant which we call a 
\textit{quantum weight} $Q$. 
More precisely, any $\mathbf{X}$-labeled unoriented tangle diagram $T$ can, 
after applying planar isotopy if necessary, be divided into pieces of the 
pictured forms:
\[\includegraphics{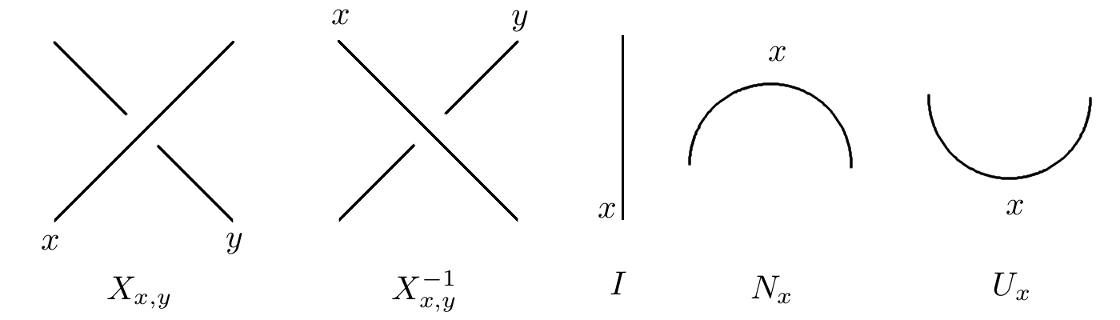}\]
Now let us fix a field $k$ and a $k$-vector space $V$. The idea is to
assign linear maps to the basic tangles so that the overall tangle
determines a linear map when we interpret 
horizontal stacking as tensor product and vertical stacking as composition 
of linear transformations.
\[\includegraphics{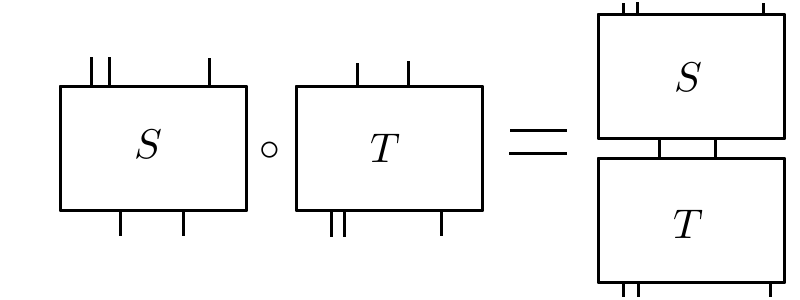}\quad \includegraphics{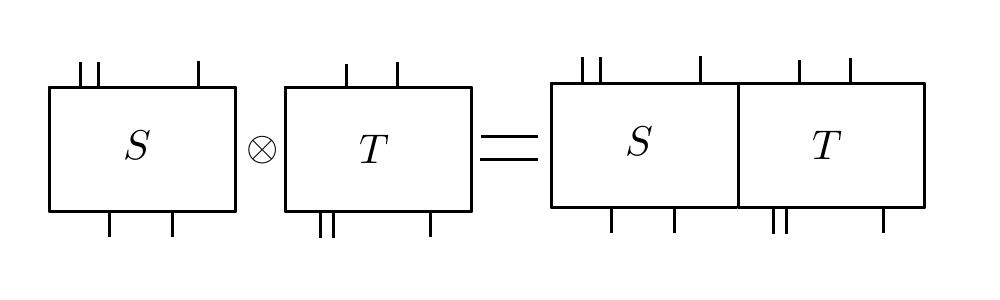}\]
A quantum weight will then be an assignment of linear transformations 
to these basic $\mathbf{X}$-labeled tangles such that equivalent 
$\mathbf{X}$-labeled framed tangles define the same linear transformation. 
This breaks down into the requirement that the assignment respects the 
\textit{$\mathbf{X}$-labeled framed tangle moves} (see \cite{SS,FY}):
\[\begin{array}{ccc}
\includegraphics{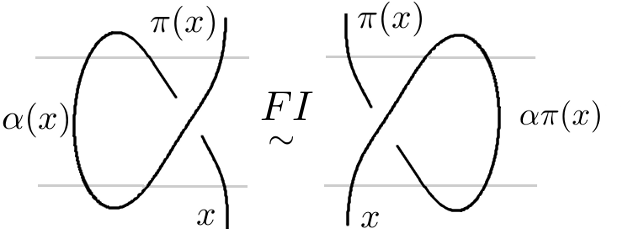} & 
\includegraphics{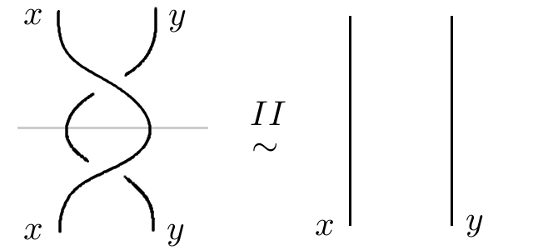} & \includegraphics{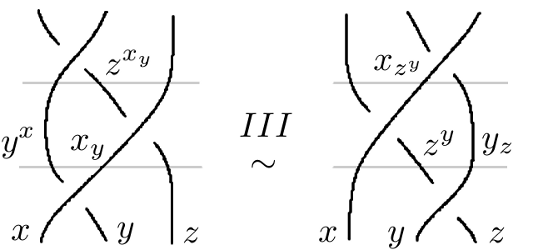} \\
\end{array}\]
\[\begin{array}{ccc}
\includegraphics{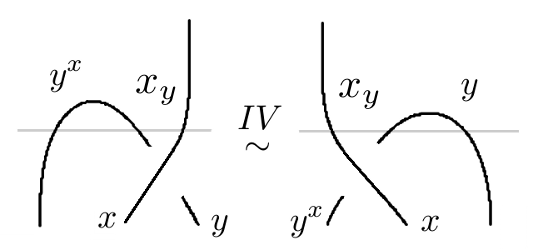} &\includegraphics{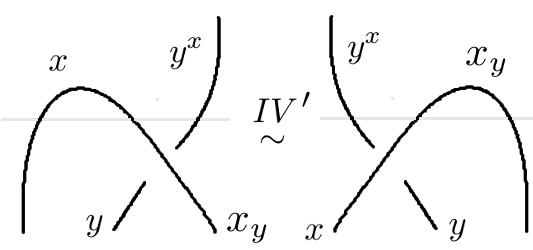} &
\includegraphics{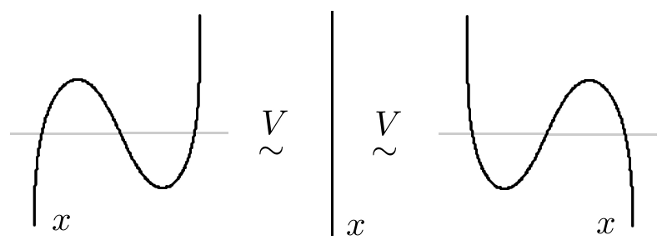} \\
\end{array}
\]
\[\begin{array}{cc}
\includegraphics{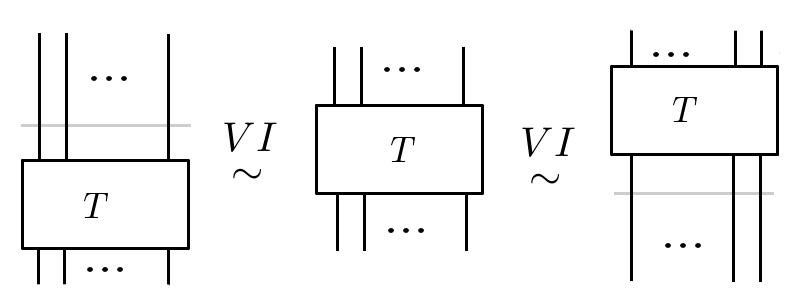} &  \includegraphics{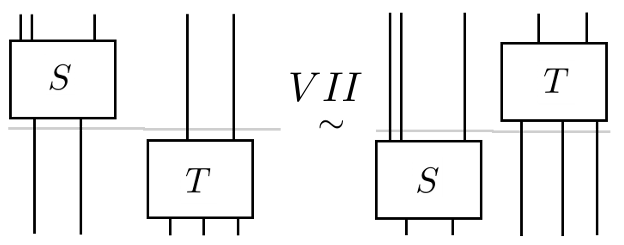}\\
\end{array}, \]
together with the requirement that the appropriate $\mathbf{X}$-labeled 
$N$-phone cord move acts as multiplication by a scalar $\delta$.
\[\raisebox{-0.35in}{\includegraphics{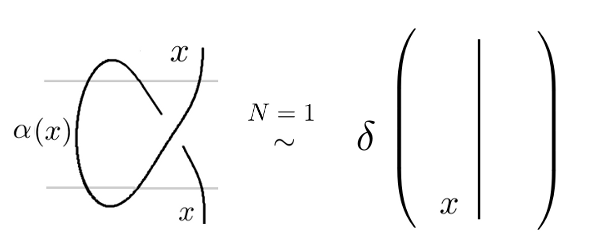}} 
\quad \mathrm{or} \quad \raisebox{-0.55in}{\includegraphics{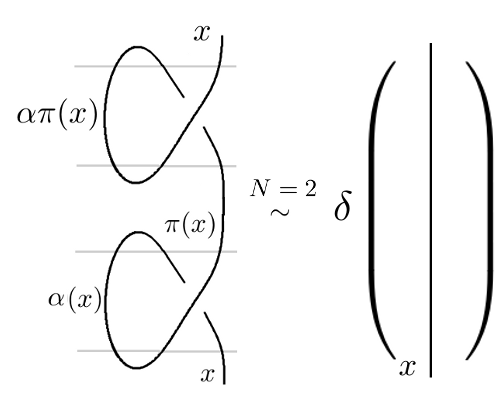}}\]

Fixing a basis for $V$, we can regard sideways stacking as Kr\"onecker 
product of matrices and vertical stacking as matrix product.
The moves $VI$ and $VII$ are automatically satisfied if we choose the identity 
transformation for the $I$-tangle by the mixed-product property of Kronecker 
product. Thus we have:

\begin{definition}
\textup{Let $\mathbf{X}$ be an involutory birack with finite rank $N$, 
$V$ a vector space over a field $k$ and $I:V\to V$ the identity transformation 
on $V$. A \textit{quantum weight} is an assignment of linear transformations}
\[X_{x,y}:V\otimes V\to V\otimes V, \quad
N_x:V\otimes V\to k, \quad
U_x:k\to V\otimes V
\]
\textup{indexed by $x,y\in \mathbf{X}$ satisfying the following conditions for
all $x,y,z\in \mathbf{X}$:}
\begin{list}{}{}
\item[(I)]
$(N_{\alpha(x)}\otimes I)(I\otimes X_{\alpha(x),x})(U_{\alpha(x)}\otimes I)
 =  (I\otimes N_{\alpha\pi(x)})(X_{x,\alpha\pi(x)}\otimes I)(I\otimes U_{\alpha\pi(x)})$
\item[(II)] $X_{x,y}$ \textup{is invertible,}
\item[(III)] $(X_{y^x,z^{x_y}}\otimes I)(I\otimes X_{x_y,z})(X_{x,y}\otimes I) 
=(I\otimes X_{x_{z^y},y_z})(X_{x,z^y}\otimes I)(I\otimes X_{y,z})$
\item[(IV)]
$(N_{y^x}\otimes I)(I\otimes X_{x,y})=(I\otimes N_y)(X_{x_y,y}^{-1}\otimes I)$
\item[(IV$'$)] $(N_x\otimes I)(I\otimes X_{x,y^x}^{-1})=(I\otimes N_{x_y})(X_{x,y}\otimes I)$
\item[(V)] $(I\otimes N_x)(U_x\otimes I)  =  I = (N_x\otimes I)(I\otimes U_x)$
\item[(VI)] $\left\{\begin{array}{ll}
(N_{\alpha(x)}\otimes I)(I\otimes X_{\alpha(x),x})(U_{\alpha(x)}\otimes I)=\delta I & N=1 \\
(N_{\alpha\pi(x)}\otimes I)(I\otimes X_{\alpha\pi(x),\pi(x)})(U_{\alpha\pi(x)}\otimes I)
(N_{\alpha(x)}\otimes I)(I\otimes X_{\alpha(x),x})(U_{\alpha(x)}\otimes I) = \delta I & N=2
\end{array}\right.$
\end{list}
\textup{for an invertible scalar $\delta$.}
\end{definition}


Given an involutory birack $\mathbf{X}=\{x_1,\dots, x_m\}$ and a $k$-vector 
space $V$ of dimension $d$ with a fixed basis, we can specify a quantum weight 
with a quadruple $(M,N,U,\delta)$ where $M$ is an $m\times m$ block matrix of 
$d^2\times d^2$ blocks, with the $(i,j)$ block of $M$ the matrix with respect
to a fixed choice of basis of the map
$X_{x_i,x_j}:V\otimes V\to V\otimes V$ associated to the crossing with labels
$x_i$ and $x_j$, $N$ is a block matrix of row vectors specifying the
maps $N_x:V\otimes V\to k$, and $U$ is a block matrix of column vectors 
specifying the maps $U_k:k\to V\otimes V$.
\[M=\left[\begin{array}{c|c|c|c}
X_{11} & X_{12} & \dots & X_{1m} \\ \hline
X_{21} & X_{22} & \dots & X_{2m} \\\hline 
\vdots & \vdots & \ddots & \vdots \\\hline
X_{m1} & X_{m2} & \dots & X_{mm} \\
\end{array}\right], \quad 
N=\left[\begin{array}{c}
N_1 \\ \hline
N_2 \\ \hline
\vdots \\ \hline
N_m
\end{array} \right],\quad \mathrm{and}\quad
U=\left[\begin{array}{c|c|c|c}
U_1 & U_2 & \dots & U_m\\ 
\end{array}\right].
\]

\begin{example}\label{ex2}\textup{
Let $B=\{1\}$, the singleton birack. Then the quadruple
\[(M,N,U,\delta)=\left(\left[\begin{array}{cccc}
A & 0 & 0 & 0 \\
0 & 0 & A^{-1} & 0 \\
0 & A^{-1} & A-A^{-3} & 0 \\
0 & 0 & 0 & A
\end{array}\right],
\quad \left[\begin{array}{cccc}
0 & A & -A^{-1} & 0 \\ 
\end{array}\right],\quad
\left[\begin{array}{c}
0 \\ -A \\ A^{-1} \\ 0 \\ 
\end{array}\right],\quad
-A^3\right)
\] defines a well-known quantum weight, the Kauffman bracket/Jones 
polynomial. Indeed, for any involutory birack $\mathbf{X}$, we can set 
$X_{x,y},N_x$ and $U_x$ equal to these $M,N$ and $U$ matrices respectively 
for all $x,y\in\mathbf{X}$ and $\delta=-A^3$ to obtain a quantum weight.
}\end{example}

\begin{example}\label{ex3}\textup{
Let $\mathbf{X}$ be the constant action birack on $\{1,2\}$ with 
$\sigma=\tau=(12)$; $\mathbf{X}$ has birack matrix
\[\left[\begin{array}{cc|cc}
2 & 2 & 2 & 2 \\
1 & 1 & 1 & 1 \\
\end{array}\right].\]
It is easy to verify that $\mathbf{X}$ is involutory (see \cite{AN}), and it is 
straightfoward (if somewhat tedious) to verify that the following quadruple
defines a quantum weight of $\mathbf{X}$ where $V=\mathbb{Q}^2$:
\[\left(\left[\begin{array}{cccc|cccc}
0 & 0 & 0 & b^{-1} & 0 & 0 & 0 & b \\
0 & a & 0 & 0 & 0 & a^{-1} & 0 & 0 \\
0 & 0 & a & 0 & 0 & 0 & a^{-1} & 0 \\
b^{-1} & 0 & 0 & 0 & b & 0 & 0 & 0 \\ \hline
0 & 0 & 0 & b & 0 & 0 & 0 & b^{-1} \\
0 & a^{-1} & 0 & 0 & 0 & a & 0 & 0 \\
0 & 0 & a^{-1} & 0 & 0 & 0 & a & 0 \\
b & 0 & 0 & 0 & b^{-1} & 0 & 0 & 0 \\
\end{array}
\right],
\left[\begin{array}{cccc}
0 & n & -n & 0 \\\hline
0 & -n & n & 0 \\
\end{array}
\right],
\left[\begin{array}{c|c}
0 & 0 \\
-n^{-1} & n^{-1} \\
n^{-1} & -n^{-1} \\
0 & 0
\end{array}
\right],-a^{-1}
\right).\]
}\end{example}

A quantum weight $Q$ translates an $\mathbf{X}$-labeling $f$ of a tangle 
diagram $T$ into a linear transformation $Q(f):V^{\otimes n}\to V^{\otimes m}$
for some $n,m\in\mathbb{Z}$. 
In the special case when $T$ is a closed tangle, i.e. a knot or link, 
$Q(f)$ is a scalar. Moreover, by construction $Q(f)$ is invariant under 
$\mathbf{X}$-labeled framed isotopy moves, and $N$-phone cord moves
change $Q(f)$ by a power of $\delta$. To obtain an unframed invariant, 
we correct for the effects of $N$-phone cord moves analogously to the 
normalization used to obtain the Jones polynomial from the Kauffman bracket:
if the blackboard framing of a component
$c_k$ of $T$ is $w_k$, write $w_k=q_kN+r_k$ where $0\le r_k<N$, and define the
\textit{normalized quantum weight} $\bar{Q(f)}$ of an $\mathbf{X}$-labeled
tangle diagram $T$ by
\[\bar{Q(f)}=\delta^{-\vec{w}}Q(f)\quad \quad \mathrm{where}\quad 
\delta^{-\vec{w}}=\prod_{k=1}^c\delta^{-q_k}.\]
Equivalently, the normalized quantum weight of a diagram is the 
quantum weight of the diagram equivalent to $T$ by framed Reidemeister and
$N$-phone cords moves whose framing is equal to its reduced value mod $N$. 
Then by construction we have

\begin{theorem} Let $\mathbf{X}$ be an involutory birack of rank $N$ and $Q$ 
a quantum weight. If $T$ and $T'$ are unoriented blackboard framed
$\mathbf{X}$-labeled tangle diagrams which are related by framed isotopy and 
$N$-phone cord moves, then $\bar{Q(L)}=\bar{Q(L')}$.
\end{theorem}

As a corollary, we can use quantum weightings to enhance the integral birack
counting invariant. More precisely, we have

\begin{definition}\textup{
Let $\mathbf{X}$ be an involutory birack, $V$ a $k$-vector space, 
$Q$ a quantum weight and $T$ an unoriented tangle of $c$ components. The 
\textit{quantum enhanced multiset invariant} of $T$ is the multiset
\[\Phi_{\mathbf{X}}^{Q,M}(L)=\{\bar{Q(f)} \ | \ 
f\in\mathrm{Hom}(IB(T,\vec{w}),\mathbf{X}), \ \vec{w}\in
\mathrm{Z}_N^c\}\]
and if $T$ is a link, the \textit{quantum enhanced polynomial invariant} of 
$T$ with respect to $\mathbf{X}$ and $Q$ is
\[\Phi_{\mathbf{X}}^Q(T)=\sum_{\vec{w}\in\mathbb{Z}_N^c} 
\left(\sum_{f\in\mathrm{Hom}(IB(T,\vec{w}),\mathbf{X})} u^{\bar{Q(f)}}\right).\]
}\end{definition}

In particular, we have
\begin{theorem}
If $\mathbf{X}$ is an involutory birack, $Q$ is a quantum weight and $T$ 
and $T'$ are ambient isotopic unoriented tangles, then 
$\Phi_{\mathbf{X}}^{Q,M}(T)=\Phi_{\mathbf{X}}^{Q,M}(T')$ and if $T$ is a link, then
$\Phi_{\mathbf{X}}^Q(T)=\Phi_{\mathbf{X}}^Q(T')$.
\end{theorem}

\begin{example}\label{ex5}
\textup{Let $\mathbf{X}$ be the birack with a single element. Then
$\Phi_{\mathbf{X}}^{Q,M}$ is a tangle functor or quantum
invariant, and indeed all such invariants can be regarded as special 
cases of $\Phi_{\mathbf{X}}^{Q,M}$ with $\mathbf{X}=\{1\}$. Alternatively, a 
quantum weight in which the maps $X_{x,y},N_x,U_x$ do not depend on the 
birack labeling satisfies 
$\Phi_{\mathbf{X}}^{Q,M}={\Phi_{\mathbf{X}}^{\mathbb{Z}}\times Q}$ where $Q$ is
a tangle functor.
}\end{example}

\begin{remark}\textup{Let $\mathbf{X}$ be a finite involutory birack of rank 
$N=1$ with $B_2(x,y)=x$, also known as a \textit{kei} 
(\begin{CJK*}{UTF8}{min}圭\end{CJK*}) or  \textit{involutory 
quandle}, and let $V=k$ be a one-dimensional vector space. Then a quantum 
weight assigns scalars $X_{x,y}$, $N_x$, and $U_x$ to each 
pair of elements or element of $\mathbf{X}$ respectively, satisfying the above 
conditions. If $N_x=U_x=\delta=1$ for all $x\in\mathbf{X}$, the function 
$\phi:\mathbf{X}\times \mathbf{X}\to k$ is then a quandle 2-cocycle and 
$\Phi_{\mathbf{X}}^Q(L)$ is the CJKLS quandle 2-cocycle 
invariant associated to $\phi$. See \cite{CJKLS} for more.
}\end{remark}

\begin{remark}\textup{
The anonymous referee of an earlier version of this paper suggested that
an alternative way to view quantum enhancements would regard them as quantum 
invariants which admit gradings by biracks. Future work will no doubt explore 
this perspective.
}\end{remark}

\begin{example}\textup{
Let $\mathbf{X}$ be the constant action birack from example \ref{ex3}. Note 
that $\mathbf{X}$-labelings of a tangle simply switch the label from
1 to 2 or 2 to 1 at every overcrossing and undercrossing point, and moreover 
$\Phi_{\mathbf{X}}^{\mathbb{Z}}(T)=2^c$ where $c$ is the number of components of $T$.
Thus, $\Phi_{\mathbf{X}}^{\mathbb{Z}}$ does not distinguish any pair of tangles with 
the same number of components. However, consider the following two-component
tangles, each with $\Phi_{\mathbf{X}}^{\mathbb{Z}}(T)=4$:
\[\includegraphics{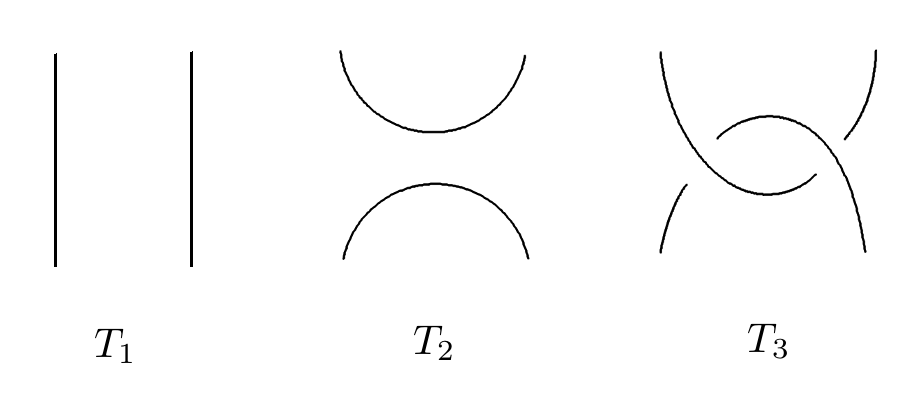}\]
The $\mathbf{X}$-labeling of $T_3$ below, for instance, contributes
\begin{eqnarray*}
\sigma\left(\raisebox{-0.3in}{\scalebox{0.7}{\includegraphics{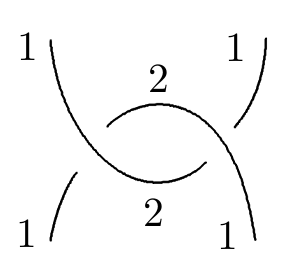}}}\right) & = & 
(I\otimes N_2\otimes I)(X_{1,2}^{-1}\otimes X_{2,1}^{-1})(I\otimes U_2\otimes I)\\
& = & \left(\left[\begin{array}{cc} 1 & 0 \\0 & 1 \end{array}\right]\otimes
\left[\begin{array}{cccc} 0 & -n & n & 0 \end{array}
\right]
\otimes\left[\begin{array}{cc} 1 & 0 \\0 & 1 \end{array}\right]
\right) \\
& & 
\left(\left[\begin{array}{cccc} 
0 & 0 & 0 & b^{-1}  \\
0 & a & 0 & 0 \\ 
0 & 0 & a & 0 \\
b^{-1} & 0 & 0 & 0
\end{array}\right]\otimes
\left[\begin{array}{cccc} 
0 & 0 & 0 & b^{-1}  \\
0 & a & 0 & 0 \\ 
0 & 0 & a & 0 \\
b^{-1} & 0 & 0 & 0
\end{array}\right]
\right) 
\left(\left[\begin{array}{cc} 1 & 0 \\0 & 1 \end{array}\right]\otimes
\left[\begin{array}{c}
0 \\ n^{-1} \\ -n^{-1} \\ 0
\end{array}\right]
\otimes\left[\begin{array}{cc} 1 & 0 \\0 & 1 \end{array}\right]\right)
\\
& = & 
\left[\begin{array}{cccc}
0 & 0 & 0 & 0 \\
0 & a^{-2} & -b^2 & 0 \\
0 & -b^2 & a^{-2} & 0 \\
0 & 0 & 0 & 0
\end{array}\right]\\
\end{eqnarray*}
to $\Phi_{\mathbf{X}}^{Q,M}(T_2)$ since $\vec{w}=(0,0)$ so $\bar{Q(f)}=Q(f)$; 
repeating for the other $\mathbf{X}$-labelings
and for the other tangles yields the invariant values listed below.
\[
\begin{array}{rcl}
\Phi_{\mathbf{X}}^{Q,M}(T_1) & = & \left\{4\times \left[\begin{array}{cccc}
1 & 0 & 0 & 0 \\
0 & 1 & 0 & 0 \\
0 & 0 & 1 & 0 \\
0 & 0 & 0 & 1
\end{array}\right]\right\},\\
\Phi_{\mathbf{X}}^{Q,M}(T_2) & = & 
\left\{2\times \left[\begin{array}{cccc}
0 & 0 & 0 & 0 \\
0 & -1 & 1 & 0 \\
0 & 1 & -1 & 0 \\
0 & 0 & 0 & 0
\end{array}\right],
2\times \left[\begin{array}{cccc}
0 & 0 & 0 & 0 \\
0 & 1 & -1 & 0 \\
0 & -1 & 1 & 0 \\
0 & 0 & 0 & 0
\end{array}\right]
\right\}, \\
\Phi_{\mathbf{X}}^{Q,M}(T_2) & = & 
\left\{2\times \left[\begin{array}{cccc}
0 & 0 & 0 & 0 \\
0 & a^{-2} & -b^2 & 0 \\
0 & -b^2 & a^{-2} & 0 \\
0 & 0 & 0 & 0
\end{array}\right],
2\times \left[\begin{array}{cccc}
0 & 0 & 0 & 0 \\
0 & -a^2 & b^{-2} & 0 \\
0 & b^{-2} & -a^2 & 0 \\
0 & 0 & 0 & 0
\end{array}\right]
\right\}. \\
\end{array}
\]
}\end{example}
Note that the invariant $\Phi_{\mathbf{X}}^Q$ determines the integral 
birack counting invariant since 
$|\Phi_{\mathbf{X}}^{M,Q}|=\Phi_{\mathbf{X}}^{\mathbb{Z}}$ and, if $T$ is a link,
evaluating $\Phi_{\mathbf{X}}^Q$ at $u=1$ yields $\Phi_{\mathbf{X}}^{\mathbb{Z}}$. As 
example \ref{ex3} shows, $\Phi_{\mathbf{X}}^Q$ is 
a stronger invariant in general than $\Phi_{\mathbf{X}}^{\mathbb{Z}}$ and thus a 
proper enhancement. 

\section{\large\textbf{Quantum Enhancements of $\mathbf{X}$-labeled Braids}}
\label{B}

Every framed oriented knot or link can be expressed as a closed braid
$\hat\beta$ for a braid $\beta=\sigma_{k_1}\dots \sigma_{k_n}$. 
Restricting our attention to closed braids yields provides 
advantages for searching for quantum weights: we need only seek matrices
$\sigma_{a,b}^{\pm 1}$ satisfying the Reidemeister III move, and we can consider
matrices with arbitrary dimension instead of being limited to perfect square
dimensions as in the tangle case. On the other hand, any pair of conjugate 
closed braids determine the same knot, so the weight matrix determined by 
an $\mathbf{X}$-labeled closed braid is not a valid signature, only its
conjugacy class. We will deal with this by taking the trace of the weight 
matrix as the signature, which is unchanged by conjugacy as well as faster
computationally than other conjugacy-class invariants such as the determinant.

\begin{definition}\textup{
Let $B_n$ be the $n$-strand braid group and let $\mathbf{X}$ be an
involutory birack. An \textit{$n$-braid weight} for $\mathbf{X}$ is an 
assignment
of an invertible matrix $\sigma_j^{x,y}$ with entries in a ring $R$ for each 
$x,y\in\mathbf{X}$ and each $j=1,\dots, n-1$ such that the 
$\mathbf{X}$-labeled braid relations 
\[\sigma_j^{x,y}\sigma_{j+1}^{x_y,z}\sigma_j^{y^x,z^{x_y}}=\sigma_{j+1}^{y,z}\sigma_j^{x,z^y}\sigma_{j+1}^{x_{z^y},y_z}\]
and
\[\sigma_j^{x,y}\sigma_k^{u,v}=\sigma_k^{u,v}\sigma_j^{x,y} \quad |j-k|<2\]
are satisfied. 
}\end{definition}

Given a braid weight for an involutory birack $X$, we can define an enhancement
of the involutory birack counting invariant for closed $n$-braids by replacing
each braid generator in an $\mathbf{X}$-labeling of $\hat\beta$ with the 
corresponding matrix and taking the 
trace as a signature. Note that an $\mathbf{X}$-labeling of a closed braid 
must have the same list of $\mathbf{X}$-labels along the top and bottom of the
braid in addition to satisfying the crossing condition at each crossing, and 
note further that fixing the braid index fixes the framing, so in this section 
we are only enhancing $\Phi^B_{\mathbb{X}}$.

\begin{definition}\textup{
Let $\beta\in B_n$, $\mathbf{X}$ be an involutory birack, and 
$W=\{\sigma_j^{x,y} \ |\ x,y\in X, 1\le j\le n-1\}$ be a braid weight for $X$. 
Then the \textit{braid weight enhancement of the basic $\mathbf{X}$-counting 
invariant} is the multiset of the traces of the matrices $\beta_f$ obtained
by replacing each braid group generator $\sigma_j$ with the appropriate matrix 
$X_j^{x,y}$ over the set of all $\mathbf{X}$-labelings $f$ of $\beta$, i.e.
\[\Phi^{M,W}_{\mathbf{X}}(\beta)=\{\trace(\beta_W)\ |\ f\in\mathrm{Hom}(IB(\beta),\mathbf{X})\}\]
with polynomial version
\[\Phi^{W}_{\mathbf{X}}(\beta)=\sum_{f\in\mathrm{Hom}(IB(\beta),\mathbf{X})} u^{\trace(\beta_W)}.\]
}\end{definition}

\begin{example}\textup{
Let $\mathbf{X}$ be the involutory birack with matrix
\[M_{\mathbf{X}}=\left[\begin{array}{cc|cc}
2 & 2 & 2 & 2 \\
1 & 1 & 1 & 1
\end{array}\right].\] 
Our computer search identified $3$-braid weights for $\mathbf{X}$ including
\[\begin{array}{|c|c|}\hline
j=1 & j=2 \\ \hline
 & \\
\begin{array}{cc}
\left[\begin{array}{cc} 0 & 1 \\ x & 0\end{array}\right] &
\left[\begin{array}{cc} 0 & 1 \\ y & 0\end{array}\right] \\  
& \\
\left[\begin{array}{cc} 0 & 1 \\ z & 0\end{array}\right] &
\left[\begin{array}{cc} 0 & 1 \\ w & 0\end{array}\right]
\end{array} &
\begin{array}{cc}
\left[\begin{array}{cc} 0 & 1 \\ w & 0\end{array}\right] &
\left[\begin{array}{cc} 0 & 1 \\ z & 0\end{array}\right] \\ 
 & \\
\left[\begin{array}{cc} 0 & 1 \\ y & 0\end{array}\right] &
\left[\begin{array}{cc} 0 & 1 \\ x & 0\end{array}\right]
\end{array} \\ & \\ \hline
\end{array}.
\]
Then for instance the closure of the braid 
$\beta=\sigma_1\sigma_1\sigma_1\sigma_2\in B_3$ is the trefoil knot with 
framing number $2$; it
has $\Phi^{M,W}_{\mathbf{X}}(\beta)=\{2y^2, 2z^2\}$ while the closure of
the braid $\beta'=\sigma_1\sigma_1^{-1}\sigma_1\sigma_2$ is an unknot
with $\Phi^{M,W}_{\mathbf{X}}(\beta)=\{2y, 2z\}$.
}\end{example}

\section{\large\textbf{Questions}}\label{Q}

We conclude with a few questions and directions for future research.

Our computations suggest that the quantum enhancement in example \ref{ex3} is
trivial on classical knots, and indeed many of the two-dimensional quantum 
enhancements we were able to find in our (far from exhaustive) computer 
search with respect to small-cardinality involutory biracks $\mathbf{X}$ 
appear to yield trivial
invariants on closed classical knots. Nevertheless, the enhancement in 
example \ref{ex3} is nontrivial on at least some tangles with nonempty 
boundary, and there are many known examples of nontrivial one-dimensional 
quantum weights (namely, CJKLS quandle 2-cocycle invariants) and 
nontrivial two-dimensional quantum weights on small cardinality biracks 
(namely, tangle functor invariants). 

Expanding on example \ref{ex5}, let us say a quantum weight is 
\textit{homogeneous} if $X_{x,y}=X_{x',y'}$, $N_x=N_{x'}$ and $U_x=U_{x'}$
for all $x,y,x',y'\in\mathbf{X}$, and that a quantum weight is 
\textit{heterogeneous}
if at least one $X_{x,y}\ne X_{x'y'}$, $N_{x}\ne N_{x'}$ or $U_x\ne U_{x'}$.
Example \ref{ex3} demonstrates the existence of heterogeneous quantum 
enhancements. For a given involutory birack $\mathbf{X}$, what is the minimal 
dimension of a vector space $V$ required for the existence of a heterogeneous
quantum weight?

Continuing in the same vein, say a quantum weight is \textit{strongly 
heterogeneous} if at least one of the $X_{i,j}$ matrices is not a classical
$R$-matrix, i.e. if $X_{i,j}$ does not satisfy the unlabeled Yang-Baxter 
equation
\[(I\times X_{i,j})(X_{i,j}\times I)(I\times X_{i,j}) =
(X_{i,j}\times I)(I\times X_{i,j})(X_{i,j}\times I).\]
If $\mathbf{X}$ is an involutory quandle, then the matrices $X_{ii}$ on the 
diagonal of $X$ must be $R$-matrices, but even in this case the off-diagonal 
matrices need not satisfy the unlabeled Yang-Baxter equation \textit{a priori},
and the quantum weights associated to a non-quandle birack might all be fail 
to be classical $R$-matrices. Such quantum weights are expected to define the 
most interesting and non-trivial quantum enhancements.

We note that in general, finding quantum weights is a difficult problem. 
Even for the smallest non-trivial biracks (those with two elements)
and the smallest non-scalar quantum weights ($\mathrm{dim}(V)=2$), the 
entries of the matrices represent up to 80 independent variables
with the axioms yielding a system of hundreds of non-linear equations. 
Computer searches have yielded some results, but better methods of finding 
quantum weights would be of great interest. Our \texttt{python} code is 
available at \texttt{www.esotericka.org}.

Finally, we note that in this paper we have considered only the simplest 
possible case, that of unoriented classical tangles. In future papers we will 
explore the oriented, virtual and twisted virtual case, each of which 
involves more complicated axioms and, we expect, richer structure.

\bigskip

\noindent
\textsc{Department of Mathematical Sciences \\
Claremont McKenna College \\
850 Columbia Ave. \\
Claremont, CA 91711}

\end{document}